\documentclass[titlepage,11pt]{article}
\usepackage{latexsym}
\usepackage{amsfonts}
\usepackage{amsmath}
\newcommand\blackslug{\hbox{\hskip 1pt \vrule width 4pt height 8pt depth 1.5pt
        \hskip 1pt}}
\newcommand\bbox{\hfill \quad \blackslug \bigbreak}

\def\ll{,\ldots,}
\def\cupcup{\cup\cdots\cup}

\title{Polynomial bounds for chromatic number \\ II. Excluding a star-forest}
\author{Alex Scott\thanks{Research supported by EPSRC grant EP/V007327/1.}\\
Mathematical Institute, University of Oxford, Oxford OX2 6GG, UK
\\
\\
Paul Seymour\thanks{Supported by AFOSR grant
A9550-19-1-0187, and by NSF grant  DMS-1800053.}\\
Princeton University, Princeton, NJ 08544
\\
\\
Sophie Spirkl\thanks{We acknowledge the support of the Natural Sciences and Engineering Research
Council of Canada (NSERC), [funding reference number RGPIN-2020-03912].
Cette recherche a \'et\'e financ\'ee par le Conseil de recherches en sciences
naturelles et en g\'enie du Canada (CRSNG), [num\'ero de r\'ef\'erence
RGPIN-2020-03912].  }\\
University of Waterloo, Waterloo, Ontario N2L3G1, Canada}

\date{}

\newtheorem{thm}{}[section]

\newcommand{\Proof}{\noindent{\bf Proof.}\ \ }

\begin{document}
\maketitle
\begin{abstract}
The Gy\'arf\'as-Sumner conjecture says that for every forest $H$, there
is a function $f_H$ such that if $G$ is $H$-free then $\chi(G)\le f_H(\omega(G))$ (where $\chi, \omega$ are the chromatic number
and the clique number of $G$). Louis Esperet conjectured that, whenever such a statement holds, 
$f_H$ can be chosen to be a polynomial. The Gy\'arf\'as-Sumner conjecture is only known to be true for a modest set of 
forests $H$, and Esperet's conjecture is known in to be true for almost no forests. For instance, it is not known when $H$ 
is a five-vertex path. Here we prove Esperet's conjecture 
when each component of $H$ is a star.
\end{abstract}

\section{Introduction}
The Gy\'arf\'as-Sumner conjecture~\cite{gyarfas, sumner} asserts:
\begin{thm}\label{GSconj}
{\bf Conjecture: } For every forest $H$, there is a function $f$
such that $\chi(G)\le f(\omega(G))$ for every $H$-free graph $G$. 
\end{thm}
(We use $\chi(G)$ and $\omega(G)$ to denote the chromatic
number and the clique number of a graph $G$, and a graph is {\em $H$-free} if it has no induced subgraph isomorphic to $H$.) 
This remains open in general, though it has been proved for some 
very restricted families of trees
(see, for example,~\cite{distantstars, gyarfasprob, gst, kierstead, kierstead2, scott, newbrooms}).

A class $\mathcal C$ of graphs is {\em $\chi$-bounded} if there is a function $f$ such that $\chi(G)\le f(\omega(G))$ 
for every graph $G$ that is an induced subgraph of a member of $\mathcal C$
 (see~\cite{survey} for a survey).  Thus the
Gy\'arf\'as-Sumner conjecture asserts that, for every forest $H$, the class of all $H$-free graphs is $\chi$-bounded. 
Esperet~\cite{esperet} conjectured that every $\chi$-bounded class is {\em polynomially $\chi$-bounded}, that is, $f$ can be chosen to be a polynomial. Neither conjecture has been settled in general. There is a survey by Schiermeyer and Randerath~\cite{Schiermeyer}
on related material.

In particular, what happens to Esperet's conjecture when we exclude a forest? For which forests $H$ can we show 
the following?
\begin{thm}\label{esperet}
{\bf Esperet's conjecture:} There is a polynomial $f_H$ such that $\chi(G)\le f_H(\omega(G))$ for every $H$-free graph $G$.
\end{thm}
Not for very many forests $H$, far fewer than the forests that we know satisfy \ref{GSconj}.
For instance, \ref{esperet} is not known when $H=P_5$, the five-vertex path. (This case is of great interest, because 
it would imply the Erd\H{o}s-Hajnal conjecture~\cite{EH0, EH, fivehole} for $P_5$, and the latter is currently the smallest open case of 
the Erd\H{o}s-Hajnal conjecture.) 

We remark that, if in \ref{esperet} we replace $\omega(G)$ by $\tau(G)$, defined to be the maximum $t$ such that $G$ contains $K_{t,t}$ as a 
subgraph, then all forests satisfy the modified \ref{esperet}. More exactly, the following is shown in~\cite{Ktt}:
\begin{thm}\label{Ktt}
For every forest $H$, there is a polynomial $f_H$ such that $\chi(G)\le f_H(\tau(G))$ for every $H$-free graph $G$.
\end{thm}

One difficulty with \ref{esperet} is that we cannot assume that $H$ is connected, or more exactly, knowing that each component of $H$  satisfies
\ref{esperet} does not seem to imply that $H$ itself satisfies \ref{esperet}. For instance, while $H=P_4$
satisfies \ref{esperet}, we do not know the same when $H$ is the disjoint union of two copies of $P_4$.

As far as we are aware, the only forests that were already known to satisfy \ref{esperet} are those of the following three 
results, and their induced subgraphs (a {\em star} is a tree in which one vertex is adjacent to all the others):
\begin{thm}\label{survey}
The forest $H$ satisfies \ref{esperet} if either:
\begin{itemize}
\item $H$ is the disjoint union of copies of $K_2$ {\rm(S. Wagon~\cite{wagon})}; or
\item $H$ is the disjoint union of $H'$ and a copy of $K_2$, and $H'$ satisfies \ref{esperet} {\rm (I. Schiermeyer~\cite{Schiermeyer2})}; or
\item $H$ is obtained from a star by subdividing one edge once {\rm (X. Liu, J. Schroeder, Z. Wang and X. Yu~\cite{liu})}.
\end{itemize}
\end{thm}
In the next paper of this series~\cite{poly3} we will show a strengthening of the third result of \ref{survey}, that is, \ref{esperet} 
is true when $H$ is a ``double star'', a tree with two internal vertices, the most general tree that does not
contain a five-vertex path.
Our main theorem in this paper is a strengthening of the second result of \ref{survey}:

\begin{thm}\label{mainthm}
If $H$ is the disjoint union of $H'$ and a star, and $H'$ satisfies \ref{esperet}, then $H$ satisfies \ref{esperet}.
\end{thm}

A {\em star-forest} is a forest in which every
component is a star.
From \ref{mainthm} and the result of~\cite{poly3}, we deduce
\begin{thm}\label{starforestthm}
If $H'$ is a double star, and $H$ is the disjoint union of $H'$ and a star-forest, 
then $H$ satisfies \ref{esperet}.
\end{thm}
As far as we know (although it seems unlikely), these might be all the forests that satisfy \ref{esperet}.

\section{The proof}
We will need the following well-known version of Ramsey's theorem:
\begin{thm}\label{Ramsey}
For $k\ge 1$ an integer, if a graph $G$ has no stable subset of size $k$, then 
$$|V(G)|\le \omega(G)^{k-1}+\omega(G)^{k-2}+\cdots+\omega(G).$$
Consequently $|V(G)|< \omega(G)^k$ if $\omega(G)>1$.
\end{thm}
\Proof The claim holds for $k\le 2$, so we assume that $k\ge 3$ and the result holds for $k-1$. Let $X$ be a clique of $G$
of cardinality $\omega(G)$, and for each $x\in X$ let $W_x$ be the set of vertices nonadjacent to $X$. From the inductive hypothesis,
$|W_x|\le \omega(G)^{k-2}+\cdots+\omega(G)$ for each $x$; but $V(G)$ is the union of the sets $W_x\cup \{x\}$ for $x\in X$,
and the result follows by adding. This proves \ref{Ramsey}.~\bbox

If $X\subseteq V(G)$, we denote the subgraph induced on $X$ by $G[X]$.
When we are working with a graph $G$ and its induced subgraphs, it is convenient to write $\chi(X)$ for $\chi(G[X])$.
Now we prove \ref{mainthm}, which we restate:
\begin{thm}\label{mainthm2}
If $H'$ satisfies \ref{esperet}, and $H$ is the disjoint union of $H'$ and a star, then $H$ satisfies \ref{esperet}.
\end{thm}
\Proof 
$H$ is the disjoint union of $H'$ and some star $S$; let $S$ have $k+1$ vertices. Since $H'$ satisfies \ref{esperet}, and 
$\chi(G)=\omega(G)$ for all graphs with $\omega(G)\le 1$, there 
exists $c'$ such that $\chi(G)\le \omega(G)^{c'}$ for every $H'$-free graph $G$.
Choose $c\ge k+2$ such that 
$$x^c-(x-1)^c\ge 1+x^{k+1}+x^{k(k+2)+c'}$$
for all $x\ge 2$ (it is easy to see that this is possible).

Let $G$ be an $H$-free graph, and write $\omega(G)=\omega$; we will show that 
$\chi(G)\le \omega^{c}$ by induction
on $\omega$. If $\omega=1$ then $\chi(G)=1$ as required, so we assume that $\omega\ge 2$.
Let $n=\omega^{k+1}$. If every vertex of $G$ has degree less than $\omega^{c}$, then the result holds as we can colour greedily, so we assume
that some vertex $v$ has degree at least $\omega^{c}$. Let $N$ be the set of all neighbours of $v$ in $G$.
Let $X_1$ be the largest clique contained in $N$; let $X_2$ be the largest clique contained in $N\setminus X_1$; and in general,
let $X_i$ be the largest clique contained in $N\setminus (X_1\cupcup X_{i-1})$. Since $|N|\ge \omega^{c}\ge n\omega$ (because $c\ge k+2$), it follows
that $X_1\ll X_{n} \ne \emptyset$. Let $X=X_1\cupcup X_{n}$, and $X_0=N\setminus X$, and $t=|X_{n}|$. Thus $1\le t\le \omega-1$
(because $\omega(G[N])<\omega$).
\\
\\
(1) {\em $\chi(N\cup \{v\})\le  t^{c}+ n\omega$.}
\\
\\
From the choice of $X_{n}$, it follows that the largest clique of $G[X_0]$ has cardinality 
at most $t<\omega$. From the inductive hypothesis, $\chi(X_0)\le t^{c}$, and since $X\cup \{v\}$ has cardinality at most
$n\omega$, it follows that $\chi(N\cup \{v\})\le  t^{c}+ n\omega$. This proves (1).

\bigskip

For each stable set $Y\subseteq X$ with $|Y|=k$, let $A_Y$ be the set of vertices in $V(G)\setminus (N\cup \{v\})$ that have
no neighbour in $Y$.  Let $A$ be the union of all the sets $A_Y$, and $B=V(G)\setminus (A\cup N\cup \{v\})$.
\\
\\
(2) {\em $\chi(A)\le (n\omega)^{k}\omega^{c'}$.}
\\
\\
For each choice of $Y$, $G[A_Y]$ is $H'$-free (because $Y\cup \{v\}$ induces a copy of $S$
with no edges to $A_Y$), and so $\chi(A_Y)\le \omega^{c'}$. Since there are at most $|X|^k\le (n\omega)^{k}$
choices of $Y$, it follows that the union $A$ of all the sets $A_Y$ has chromatic number at most  
$(n\omega)^{k}\omega^{c'}$. This proves (2).
\\
\\
(3) {\em For each $b\in B$, $b$ has fewer than $\omega^k$ non-neighbours in $X$.}
\\
\\
Let $Z$ be the set of vertices in $X$ nonadjacent to $b$. Since $b\notin A$, $G[Z]$
has no stable set of cardinality $k$; and since it also has no clique of cardinality $\omega$,
\ref{Ramsey} implies that $|Z|\le (\omega-1)^k<\omega^k$. This proves (3).
\\
\\
(4) {\em $\chi(B)\le (\omega-t)^c$.}
\\
\\
Suppose that $C\subseteq B$ is a clique with $|C|=\omega-t+1$. For each $c\in C$, (3) implies that
$c$ has a nonneighbour in fewer than $\omega^k$ of the cliques $X_1\ll X_n$; and so fewer than $(\omega-t+1)\omega^k$
of the cliques $X_\ll X_n$ contain a vertex with a non-neighbour in $C$. Since $(\omega-t+1)\omega^k\le \omega^{k+1}=n$,
there exists $i\in \{1\ll n\}$ such that every vertex in $X_i$ is adjacent to every vertex of $C$, and so $C\cup X_i$
is a clique. Since $|X_i|\ge |X_n|=t$, it follows that $|C\cup X_i|>\omega$, a contradiction.
Thus there is no such clique $C$, and so $\omega(G[B])\le \omega-t$; and from the inductive hypothesis (since $t>0$)
it follows that $\chi(B)\le (\omega-t)^c$. This proves (4).

\bigskip
From (1), (2), (4) we deduce that
$$\chi(G)\le \chi(N\cup \{v\})+\chi(A)+\chi(B)\le t^c+ n\omega+(n\omega)^{k}\omega^{c'}+(\omega-t)^{c}.$$
Since $1\le t\le \omega-1$ and $c\ge 1$, it follows that $t^{c}+(\omega-t)^{c}\le 1 + (\omega-1)^{c}$,
and so 
$$\chi(G)\le  1+ \omega^{k+1}+(n\omega)^{k}\omega^{c'}+(\omega-1)^{c}\le \omega^c$$
from the choice of $c$ and since $\omega\ge 2$.
This proves \ref{mainthm}.~\bbox


\begin{thebibliography}{99}
\bibitem{distantstars} 
M. Chudnovsky, A. Scott and P. Seymour, 
``Induced subgraphs of graphs with large chromatic number.  XII. Distant stars'', 
{\em J. Graph Theory} {\bf 92} (2019), 237--254,
{\tt arXiv:1711.08612}.
\bibitem{fivehole} M. Chudnovsky,  A. Scott, P. Seymour and S. Spirkl,  ``Erd\H{o}s-Hajnal for graphs with no five-hole'', submitted for publication, {\tt arXiv:2102.04994}.
\bibitem{EH0} P. Erd\H{o}s and A. Hajnal, ``On spanned subgraphs of graphs'',
{\em Graphentheorie und Ihre Anwendungen} (Oberhof, 1977). 
\bibitem{EH}  P. Erd\H{o}s and A. Hajnal, ``Ramsey-type theorems'',
{\em  Discrete Applied Math.} {\bf 25} (1989), 37--52.
\bibitem{esperet} L. Esperet,
{\em Graph Colorings, Flows and Perfect Matchings},
Habilitation thesis, Universit\'e Grenoble Alpes (2017), 24, https://tel.archives-ouvertes.fr/tel-01850463/document.
\bibitem{gyarfas}
A. Gy\'arf\'as, ``On Ramsey covering-numbers'',
in {\em Infinite and Finite Sets, Vol. II} (Colloq., Keszthely, 1973), {\em Coll. Math. Soc. J\'anos Bolyai} {\bf 10}, 801--816.
\bibitem{gyarfasprob}
A. Gy\'arf\'as, 
``Problems from the world surrounding perfect graphs'', 
{\em Proceedings of the International Conference on Combinatorial Analysis and its Applications},  (Pokrzywna, 1985),
{\em Zastos. Mat.}  {\bf 19} (1987), 413--441.
\bibitem{gst}
A. Gy\'arf\'as, E. Szemer\'edi and Zs. Tuza,
``Induced subtrees in graphs of large chromatic number'',
{\em Discrete Math.} {\bf 30} (1980), 235--344.
\bibitem{kierstead}
H. A. Kierstead and S.G. Penrice,
``Radius two trees specify $\chi$-bounded classes'',
{\em J. Graph Theory} {\bf 18} (1994), 119-–129.
\bibitem{rodl}
H. A. Kierstead and V. R\"odl, ``Applications of hypergraph coloring to coloring graphs
not inducing certain trees'', 
{\em Discrete Math.} {\bf 150} (1996), 187--193.
\bibitem{kierstead2} 
H. A. Kierstead and  Y. Zhu, 
``Radius three trees in graphs with large chromatic number'',
{\em SIAM J. Disc. Math.} {\bf 17} (2004), 571--581.
\bibitem{liu} X. Liu, J. Schroeder, Z. Wang and X. Yu, ``Polynomial $\chi$-binding functions for $t$-broom-free graphs'',
{\tt arXiv:2106.08871}.
\bibitem{scott}
A. Scott, 
``Induced trees in graphs of large chromatic number'',
{\em J. Graph Theory} {\bf 24} (1997), 297--311.
\bibitem{newbrooms} 
A. Scott and P. Seymour, 
``Induced subgraphs of graphs with large chromatic number.  XIII. New brooms'', 
{\em European J. Combinatorics} {\bf 84} (2020), 103024, {\tt arXiv:1807.03768}.
\bibitem{survey}
A. Scott and P. Seymour,
``A survey of $\chi$-boundedness'', {\em J. Graph Theory} {\bf 95} (2020), 473--504, {\tt arXiv:1812.07500}.
\bibitem{Ktt}
A. Scott, P. Seymour and S. Spirkl, ``Polynomial bounds for chromatic number. I.
Excluding a biclique and an induced tree'', submitted for publication,
{\tt arXiv:2104.07927}.
\bibitem{poly3} A. Scott, P. Seymour and S. Spirkl, ``Polynomial bounds for chromatic number. III. Excluding a double star'', in preparation.
\bibitem{Schiermeyer2}
I. Schiermeyer, ``On the chromatic number of ($P_5,$ windmill)-free graphs'', {\em Opuscula Math.} {\bf 37} (2017), 
609--615.
\bibitem{Schiermeyer}
I. Schiermeyer and B. Randerath, ``Polynomial $\chi$-binding functions and forbidden induced subgraphs: a survey'',
{\em Graphs and Combinatorics} {\bf 35} (2019), 1--31.
\bibitem{sumner}
D. P. Sumner,
``Subtrees of a graph and chromatic number'', in
{\em The Theory and Applications of Graphs}, (G. Chartrand, ed.),
John Wiley \& Sons, New York (1981), 557--576.
\bibitem{wagon} S. Wagon, ``A bound on the chromatic number of graphs without certain induced subgraphs'',
{\em J. Combinatorial Theory, Ser. B}, {\bf 29} (1980), 345--346.


\end{thebibliography}
\end{document}